\documentclass{svjour3}
\usepackage[utf8]{inputenc}
\usepackage{bm}
\usepackage{amssymb}
\usepackage{graphicx}
\usepackage{epstopdf}
\usepackage{amsfonts}
\usepackage[figuresright]{rotating}  
\usepackage{amssymb}
\usepackage{amsmath}
\usepackage{url}
\usepackage{mathrsfs} 
\usepackage{comment}
\usepackage{graphicx}

\usepackage[dvipsnames, usenames]{xcolor}

\title{A GPU-accelerated mixed-precision WENO method for extremal black hole and gravitational wave physics computations}
\author{Scott E. Field \and Sigal Gottlieb \and Zachary J. Grant \and  Leah F. Isherwood \and Gaurav Khanna}
\institute{S.E. Field \and S. Gottlieb \and L.F. Isherwood \and G. Khanna 
\email{gkhanna@umassd.edu}
\at University of Massachusetts Dartmouth, North Dartmouth, MA 02747
\and
Z. J. Grant 
\at 
Department of Computational and Applied Mathematics, Oak Ridge National Laboratory, Oak Ridge TN 37830
 }

\begin{document}
\maketitle

\begin{abstract}
We develop and use a novel mixed-precision weighted
essentially non-oscillatory (WENO) method 
for solving the Teukolsky equation, which arises when modeling
perturbations of Kerr black holes.
We show that WENO methods outperform higher-order finite-difference methods,
standard in the discretization of the Teukolsky equation, due to the 
need to add dissipation for stability purposes in the latter. 
In particular, as the WENO scheme uses no additional dissipation it is 
well-suited for scenarios requiring long-time evolution such 
as the study of Price tails and gravitational wave emission from extreme mass ratio binaries.
In the mixed-precision approach, the expensive computation of the 
WENO {\em weights} is performed in reduced floating-point precision that results in a significant 
speedup factor of $\approx 3.3$. 
In addition, we use state-of-the-art Nvidia 
general-purpose graphics processing units
and cluster parallelism to further accelerate the WENO computations. 
Our optimized WENO solver can be used to quickly generate accurate results
of significance in the field of black hole and gravitational wave physics. 
We apply our solver to study the behavior of the
Aretakis charge -- a conserved quantity, that if 
detected by a gravitational wave observatory like LIGO/Virgo would prove
the existence of extremal black holes. 
\end{abstract}

\section{Introduction}
The field of computational relativity is undergoing a transformative renaissance of sorts, due to multiple new discoveries such as the direct detection of gravitational waves from a black hole binary system by LIGO in 2015~\cite{Abbott:2016blz}, which was a awarded a Nobel Prize in 2017, the first-ever image of a black hole event horizon by the Event Horizon Telescope in 2019~\cite{akiyama2019first}, and others. Accurate computational models are critical for the success of such efforts since high-fidelity models enable improved data analysis efforts 
required for parameter estimation, tests of general relativity, and studying the formation 
channels of compact binaries~\cite{gwastro-PopulationReconstruct-Parametric-Wysocki2018,TheLIGOScientific:2016src,LIGO-O2-Catalog}.

Black hole perturbation theory is the standard framework for studying gravitational
phenomena, such as gravitational waves and linear stability of
black hole solutions. The theory of such perturbations is well developed, and
starts with pioneering investigations by Regge and Wheeler~\cite{PhysRev.108.1063} for 
Schwarzschild (nonspinning) black holes. The theory was later extended by Teukolsky
to handle perturbations of the Kerr metric~\cite{teukolsky1973perturbations},
which is an exact solution of the Einstein equation. The Kerr solution
is parameterized by the black hole mass, $M$, and the hole's spin, $a$,
and small perturbations about this solution obey the Teukolsky master equation Eqn.~\eqref{teuk0}.

For many realistic problems, we require highly-accurate numerical solutions to the
Teukolsky equation over extremely long simulation times. Yet most state-of-the-art 
solvers are based on a finite-difference numerical evolution scheme 
stabilized with numerical
dissipation~\cite{sundararajan2007towards,sundararajan2008towards,sundararajan2010binary,zenginouglu2011null}, which,
as we will show, can significantly degrade the quality of the numerical solution.

In this paper, we develop and use a novel mixed-precision weighted
essentially non-oscillatory (WENO) method for solving the Teukolsky equation. 
Our work is one of the first applications of WENO methods to the modeling of gravitational waves, 
and the first WENO-based Teukolsky solver. To showcase the improvements of this new solver,
we compare a
%few numerical approaches in this work, namely a 
standard sixth-order centered differencing approach (with an 
eighth-order Kreiss-Oliger dissipation operator~\cite{Kreiss73} for numerical stability) and third-order and 
fifth-order weighted essentially non-oscillatory (WENO) schemes. 
We show that the fifth order WENO of Jiang and Shu \cite{jiang1996efficient} performs well in the sense of stability and preservation 
of accuracy over long evolutions.

While our new numerical solver can be used to study a variety of interesting questions in gravitational physics,
in this paper, we focus our attention on a specific computational problem in the area of black hole physics. 
The physical scenario one may imagine is that of an isolated rotating, Kerr black hole, spinning at
the maximum possible rate, i.e. an extremal Kerr black hole. If such a black hole is perturbed because of any nearby 
material, or scalar or gravitational wave, 
%rigorous results from mathematical physicists such as Aretakis
analytic results suggest that the hole will undergo a rather unexpected evolution that would result in a unique signal which, if 
measured, would be a ``signature''  of an extremal black hole~\cite{Aretakis_2013,Angelopoulos:2018yvt}. In particular, Aretakis  
derived a mathematical quantity, often called ``charge'' or ``horizon hair'', that is conserved (stays constant in time) 
in such a system and this quantity is measurable from a far away distance. Black holes that are not extremally spinning
are unable to preserve the constancy of such a quantity~\cite{Burko:2019qqx}. 

Aretakis had derived those key results in the context of scalar waves, that are simpler to work with. However, 
the key ideas of Aretakis have been extended to the astrophysically realistic context of gravitational waves 
using a computational approach by Khanna and collaborators~\cite{burko2020scalar}. 
In this paper, we focus entirely on the realistic gravitational waves scenario and discuss the 
computational challenges involved with the relevant simulations and how to solve them. 
The main challenge that is involved is the accuracy of computed numerical solution and 
being able maintain that accuracy in long evolutions. As a simple example, one may compute the Aretakis conserved 
charge numerically and study how long the numerical simulation is able to maintain its constancy (within some 
tolerance). That is a key test that we use in this work. There are also additional similar tests based on the 
expected temporal behavior of the spatial derivatives of the solution that we conduct in this work.  
The numerical results from our WENO simulations retain the predicted behavior of the
Aretakis conserved charge in fidelity with rigorous mathematical results.  
The simulations are computationally costly, and take many months of computational time. To alleviate this 
problem, we implement a GPU-accelerated version of this code on a GPU-cluster, IBM/MIT's {\em Satori} a 64-node 
GPGPU system with 256 Nvidia V100 GPUs. Furthermore, we introduce and use a mixed precision WENO algorithm which 
speeds up the computation by a factor of between 3 -- 4 without meaningfully reducing its accuracy.

%to simulate field perturbations in the space-time of Kerr black holes.  
%Computational relativity is a subfield of the broad area of gravitational physics that focuses on numerical 
%simulations of physical systems such as black holes, gravitational waves, neutron stars, cosmological models 
%and so on. This field 

\section{Numerical Solution of the Teukolsky Equation}
In this section we briefly describe the Teukolsky equation~\eqref{teuk0}
and the context in which it arises.
We provide a description of the coordinate-systems used, the relevant evolution equations,
and a fully first-order reduction of the Teukolsky equation~\eqref{eq:evln}
which is subsequently discretized using the WENO method presented in Sec.~\ref{sec:computational_approach}. 

The Kerr metric, which describes a rotating black hole, 
is both an astrophysically and theoretically important exact solution of the Einstein equation.
Studying small perturbations of the Kerr metric is an especially important line 
of inquiry and can be used, for example, to understand the 
behavior of gravitational waves in the spacetime of a rotating black hole. 
The coordinate system that is typically used to describe 
the spacetime of black holes is the Boyer-Lindquist system, $(t,r,\theta,\varphi)$, 
that has many similarities with spherical coordinates. However, 
because they suffer from a coordinate-singularity at the horizon locations,
Boyer-Lindquist coordinates are not well for numerical computations. We 
instead make use of a better suited coordinate system known as ingoing Kerr coordinates. 
These are a Kerr spacetime generalization of the better-known Eddington-Finkelstein 
coordinates that are able to smoothly ``penetrate'' the horizon of a black hole. 
In the following subsection, we review the relationship between these 
different coordinate systems and emphasize some of their important aspects.  

The main evolution equation of interest in this work is the Teukolsky master 
equation that describes scalar, vector and tensor field perturbations in the 
spacetime of a Kerr black hole ~\cite{teukolsky1973perturbations} to linear order. We numerically solve this 
equation for the gravitational waves case using a compactified form of the ingoing Kerr 
coordinates. Using  hyperboloidal compactification allows us to {\em directly} sample 
the behavior of fields throughout the spacetime, including even (null) infinity. 
One important aspect of this work is that we must evolve the fields for a 
long duration because we are interested in the late-time, power-law decay behavior of 
these fields. This behavior typically appears {\em after} the quasi-normal modes of the 
system have exponentially decayed enough to become subdominant. This poses certain 
challenges that are explained in some detail in the following sections.  

The following subsections offer additional details including the main expressions 
for the quantities involved and also our computational methodology. 

\subsection{Teukolsky Equation}
The Teukolsky master equation describes scalar, vector and tensor field perturbations in the space-time of
Kerr black holes~\cite{teukolsky1973perturbations}. In Boyer-Lindquist coordinates, this equation takes the form
\begin{eqnarray}
\label{teuk0}
&&
-\left[\frac{(r^2 + a^2)^2 }{\Delta}-a^2\sin^2\theta\right]
         \partial_{tt}\Psi
-\frac{4 M a r}{\Delta}
         \partial_{t\phi}\Psi \nonumber \\
&&- 2s\left[r-\frac{M(r^2-a^2)}{\Delta}+ia\cos\theta\right]
         \partial_t\Psi\nonumber\\  
&&
+\,\Delta^{-s}\partial_r\left(\Delta^{s+1}\partial_r\Psi\right)
+\frac{1}{\sin\theta}\partial_\theta
\left(\sin\theta\partial_\theta\Psi\right)+\nonumber\\
&& \left[\frac{1}{\sin^2\theta}-\frac{a^2}{\Delta}\right] 
\partial_{\phi\phi}\Psi +\, 2s \left[\frac{a (r-M)}{\Delta} 
+ \frac{i \cos\theta}{\sin^2\theta}\right] \partial_\phi\Psi  \nonumber\\
&&- \left(s^2 \cot^2\theta - s \right) \Psi = 0  ,
\end{eqnarray}
where $M$ is the mass of the black hole, $a$ its angular momentum per unit mass, $\Delta = r^2 - 2 M r + a^2$ and 
$s$ is the ``spin weight'' of the field. The $s = \pm 2$ versions of these equations describe the radiative degrees 
of freedom of the gravitational field, and thus are the equations of interest here. As mentioned previously, this equation 
is an example of linear, hyperbolic, homogeneous PDEs which are quite common in several areas of science and 
engineering. Finite-difference methods are by far the most common numerical schemes developed for the Teukolsky Equation.

\subsection{Teukolsky Equation in Ingoing-Kerr Coordinates}
We begin with an expression of the usual Boyer-Lindquist coordinate version of 
the Kerr spacetime metric and the associated Teukolsky equation~\cite{teukolsky1973perturbations}. 
It is clear that the equation exhibits pathological behavior at the horizon 
locations, i.e. when $\Delta=0$. Note that this coordinate singularity can be 
easily removed by a suitable change of coordinates. 

To remove the coordinate singularity at the horizon locations determined by $\Delta=0$, 
we consider the above equations in the ingoing Kerr coordinate system 
$({\tilde t},r,\theta, {\tilde \varphi})$, also called ``horizon penetrating'' 
coordinates. 
These coordinates are related to the Boyer-Lindquist coordinates through 
the transformations ${\tilde \varphi}=\varphi+\int a\Delta^{-1}\,dr$ and 
${\tilde t}=t-r+r_*$, where the ``tortoise'' radial coordinate 
$r_*=\int(r^2+a^2)\Delta^{-1}\,dr$. This system does not suffer from any 
pathologies at the horizon locations and is therefore well-suited for analyzing 
fields both in the exterior and interior spacetimes of a rotating black hole. 
It is also useful to note that at the horizon location at late times, 
the $\tilde t$ variable is essentially the null variable $v = t + r_*$.

The last ingredient that goes into the setup of our coordinate system is 
hyperboloidal compactification as developed by Zengino\v{g}lu~\cite{Zenginoglu:2009hd,Zenginoglu:2011jz,Zenginoglu:2010cq,Zenginoglu:2009ey,Zenginoglu:2008uc,Zenginoglu:2008pw,Zenginoglu:2008wc,Zenginoglu:2007jw,Zenginoglu:2012us}. 
We define a compactified coordinate system $(\tau,\rho,\theta,{\tilde\varphi})$ 
by
\begin{equation}
\tau = {\tilde{t}} - {r}^2/({r}+S) + 4 \ln [S/({r}+S)]
\end{equation}
and
\begin{equation}
\rho ={r}/[1+{r}/S]
\end{equation}
where a free parameter $S$ controls both the domain and also the foliation.  
Note that $\rho\in [0,S)$ maps $r\in [0,\infty )$ and is therefore a one-to-one 
compactifying coordinate. A Penrose diagram of the slices defined by these 
coordinates in the Kerr spacetime context can be found in Ref.~\cite{Burko:2016uvr}. 
We do not show the final form of the Teukolsky master equation in these 
compactified coordinates because of the lengthy nature of the expression and 
the fact that it is not particularly illuminating. We simply refer to it in the 
symbolic form,
%\scott{[In the previous equation we had $\Psi$. Is this the same as $\psi$ below or are these related through a formula?]}
%\gk{[Good catch! Fixed.]}
\begin{eqnarray} \label{eq:near_final} 
{A}^{\tau \tau}\partial_\tau^2 \Psi &+& {A}^{\tau \rho}\partial_\tau\partial_{\rho}\Psi + {A}^{\rho \rho}\partial_{\rho}^2 \Psi + {A}^{\theta\theta} \partial_\theta^2 \Psi \nonumber \\
&+& {B}^\tau\partial_\tau\Psi + {B}^{\rho} \partial_{\rho} \Psi + {B}^\theta \partial_\theta \Psi + {C} \Psi = 0, 
\end{eqnarray}
and refer the reader to 
%the recent relevant research literature 
Refs.~\cite{Zenginoglu:2011zz,Harms:2014dqa} wherein additional details may be found.   

The computational grid is defined as a uniform grid over the compactified $\rho$ 
coordinate. As pointed out earlier, this allows us to access null infinity 
directly on the computational grid ($\rho = S$ maps to infinity).  Moreover, the 
compactification offers a solution to the  ``outer boundary problem'' 
in numerical relativity. Typical boundary conditions used in the research community 
lead to spurious wave reflections from the edge of the computational grid. However, 
with the approach of hyperboloidal compactification, one is able to extend the 
computational domain to infinity, making it possible to completely eliminate any such 
reflections~\cite{Zenginoglu:2011zz}\footnote{Note that a 
different approach towards compactification was used therein -- a hyperboloidal compactified 
{\em layer} was attached to the outer part of a Boyer-Lindquist coordinates based computational 
grid.}. In addition, the compactification allows us to employ a 
very dense computational grid (typically, $S\sim 20$) which results in highly accurate 
numerical results. Those details are provided in the next subsection.

\subsection{Additional Implementation Remarks}
After performing the transformations presented in the previous sections, we rewrite the vacuum equation in the form
\begin{eqnarray} \label{eq:final} 
\partial_\tau^2 \Psi &=& \tilde{A}^{\tau \rho}\partial_\tau\partial_{\rho}\Psi + \tilde{A}^{\rho \rho}\partial_{\rho}^2 \Psi + \tilde{A}^{\theta\theta} \partial_\theta^2 \Psi \nonumber \\
&+& \tilde{B}^\tau\partial_\tau\Psi + \tilde{B}^{\rho} \partial_{\rho} \Psi + \tilde{B}^\theta \partial_\theta \Psi + \tilde{C} \Psi, 
\end{eqnarray}
where the coefficients with a tilde are obtained by dividing the coefficients of Teukolsky equation \eqref{eq:near_final} by $-A^{\tau\tau}$. We put the 
equation \eqref{eq:final} in first-order form (in $\rho$ and $\tau$) by defining 
two new field variables
\begin{eqnarray}
\pi &\equiv& \partial_{\tau}{\Psi} + b \, \partial_{\rho}\Psi \; , \\
b & \equiv & - (\tilde{A}^{\tau \rho} + \sqrt{ (\tilde{A}^{\tau \rho})^{2} + 4 \tilde{A}^{\rho \rho}})/2 \; .
\end{eqnarray}
We chose
%to use this
these auxiliary variables because it has been discovered through extensive experimentation that the resulting first-order form,
\begin{eqnarray}
\label{eq:evln}
\partial_{\tau} \mbox{\boldmath{$u$}} + \mbox{\boldmath{$M$}} \partial_{\rho}\mbox{\boldmath{$u$}} 
+ \mbox{\boldmath{$Lu$}} + \mbox{\boldmath{$Au$}} = 0,
\end{eqnarray}
is ideally suited for long stable evolutions.
%first-order form of the equation for our numerical evolutions, because it has been discovered through extensive experimentation that such a form is ideally suited for long stable evolutions.
%Then the equation takes the form 
%\begin{eqnarray}
%\label{eq:evln}
%\partial_{\tau} \mbox{\boldmath{$u$}} + \mbox{\boldmath{$M$}} \partial_{\rho}\mbox{\boldmath{$u$}} 
%+ \mbox{\boldmath{$Lu$}} + \mbox{\boldmath{$Au$}} = 0,
%\end{eqnarray}
%where 
Here
\begin{equation}
\mbox{\boldmath{$u$}}\equiv\{\Psi_R,\Psi_I,\pi_R,\pi_I\}
\end{equation}
is the solution vector and the subscripts $R$ and $I$ refer to the real and imaginary parts respectively. The matrices {\boldmath{$M$}}, {\boldmath{$A$}} and {\boldmath{$L$}} are obtained from \eqref{eq:final} as in \cite{sundararajan2007towards,sundararajan2008towards}. Here it will suffice to simply indicate the final form taken by these matrices: 
\begin{equation}
\mbox{\boldmath{$M$}} \equiv \left(\begin{matrix}
                    b  &   0   &  0     &  0 \cr
                    0  &   b   &  0     &  0 \cr
                    m_{31}  &   m_{32}  & -b  &  0 \cr
                    -m_{32}  &   m_{31} &  0  & -b \cr
                \end{matrix}\right) \; ,
\label{m_matrix}
\end{equation}
\begin{equation}
\mbox{\boldmath{$A$}} \equiv \left(\begin{matrix}
                    0  &   0   &  -1  &  0 \cr
                    0  &   0   &  0  &  -1 \cr
                    a_{31} & a_{32} & a_{33} & a_{34} \cr
                    -a_{32} & a_{31} & -a_{34} & a_{33} \cr
                \end{matrix}\right) \; ,
\label{a_matrix}
\end{equation}
and
\begin{equation}
 \mbox{\boldmath{$L$}} \equiv \left(\begin{matrix}
                    0  &   0   &  0  &  0 \cr
                    0  &   0   &  0  &  0 \cr
                    l_{31}  &   0   &  0  &  0 \cr
                    0  &   l_{31}   &  0  &  0 \cr
                \end{matrix}\right)\;.
\label{l_matrix}
\end{equation}
The angular derivatives are encoded in {\boldmath{$L$}}.

The main advantage of casting the equation \eqref{eq:final} into the form \eqref{eq:evln} is that the system has advantageous properties in the variable $\rho$. The matrix {\boldmath{$M$}} has a complete set of linearly independent eigenvectors with real eigenvalues. This is not a rigorous statement on the hyperbolicity of the system because the matrix {\boldmath{$L$}} contains second-order angular derivatives. Nevertheless, experiments show that the system is numerically well-behaved. 

%To develop a numerical solver scheme for this equation, we write \eqref{eq:evln} as 
%\begin{equation}
%\partial_{\tau} \mbox{\boldmath{$u$}} + \mbox{\boldmath{$D$}}
%\partial_{\rho} \mbox{\boldmath{$u$}}
%=  \mbox{\boldmath{$S$}}\; , 
%\label{new_teu2}
%\end{equation}
%where
%\begin{equation}
% \mbox{\boldmath{$D$}} \equiv \left(\begin{matrix}
%                    b &   0   &  0  &  0 \cr
%                    0  &   b   &  0  &  0 \cr
%                    0  &   0   &  -b  &  0 \cr
%                    0  &   0   &  0  &  -b \cr
%                \end{matrix}\right),
%\label{d_matrix}
%\end{equation}
%\begin{equation}
%\mbox{\boldmath{$S$}} = -(\mbox{\boldmath{$M$}} - \mbox{\boldmath{$D$}})
%\partial_{\rho}\mbox{\boldmath{$u$}}
%- \mbox{\boldmath{$L$}}\mbox{\boldmath{$u$}} 
%- \mbox{\boldmath{$A$}}\mbox{\boldmath{$u$}}  .
%\end{equation}
%
\section{Computational Approach} \label{sec:computational_approach}
The numerical approach used to solve the first-order Teukolsky equation 
\eqref{eq:evln} in the compactified ingoing Kerr coordinate
system is very similar to the one presented in our earlier 
work~\cite{Zenginoglu:2011zz}. We simply outline the main steps here and refer the reader to 
that reference for additional details. We begin by taking advantage of Kerr 
spacetime's axisymmetry and separating out the $\tilde\varphi$ dependence of the 
system using an $\exp(im{\tilde\varphi})$ form for the gravitational field $\Psi$. This 
transforms the original (3+1)D equation into a system of (2+1)D equations. In this 
work we restrict ourselves to axisymmetric fields only, and therefore we set $m=0$ 
throughout. Next, we cast the equations into a first-order hyperbolic partial differential 
equation form, by defining a new ``momentum'' field that is related to the derivative 
of the field $\Psi$. 

In the next subsections we will describe the numerical methods used to discretize these 
equations. We implement a fifth-order WENO finite-difference  and a third-order WENO 
finite-difference scheme, and compare it with a sixth-order finite-difference scheme with 
numerical dissipation required for stability. The time-stepping method
used is the  third-order Shu-Osher explicit Runge--Kutta scheme. 

It is worth commenting on the fact that numerical computations are rather  
challenging in the context of studying the late-time tails. As remarked before, 
these computations must be long duration because the observed field initially exhibits 
an exponentially decaying oscillatory behavior known as quasi-normal ringing. Only much 
later, once the exponential decay has made these modes subdominant, does the field   
transition over to a power-law tail. Moreover, there are often intermediate tails~\cite{Zenginoglu:2012us}, 
that do not necessarily have the true late-time asymptotic rates that we are interested 
in here. These intermediate tails decay faster than the asymptotic rate, but may 
have dominant amplitudes for a period of time. We must evolve longer than these 
intermediate tails last in order to obtain the tail solution with the true asymptotic 
decay rate. 

In addition, each of the field's spherical harmonic multipoles $Y_{\ell m}$ has its 
own decay rate (that is proportional to $\ell$). Thus, at late times we obtain 
numerical data in which different multipoles may have widely ranging amplitudes 
(typically 15 -- 20 orders of magnitude apart!). It is thus important for the numerical 
scheme to have high grid density in order to reduce the {\em truncation} errors 
to very low levels. In addition, due to the very large range of amplitudes involved, 
these computations also require high-precision floating-point numerics that allow us 
to reduce {\em round-off} error that can otherwise easily overwhelm the fast decaying 
multipoles. In particular, we satisfy this requirement by using {\em quad}-precision 
numerics (128-bit or $\sim$30 decimal digits). This keeps the round-off error in our 
computations at acceptably low levels. 

Finally, to complete these long duration, high-accuracy and high-precision computations 
in a reasonable time-frame we make extensive use GPGPU-based parallel computing. For 
additional details on implementation of such intensive computations on a parallel GPU 
architecture, we refer the reader to our earlier work on the subject~\cite{khanna2013high}.

\subsection{Spatial Discretization}
We compare two approaches in this work: a sixth-order finite difference method that is stabilized by an eighth-order Kreiss-Oliger 
dissipation operator~\cite{Kreiss73} and the weighed essentially non-oscillatory methods of Jiang and Shu \cite{jiang1996efficient}.
%\scott{[Q: Below we use x as the independent variable. Should this be $\rho$? ]} \gk{[Sigal .. your call -- I think I fixed it,please check]}
%\scott{[Q: Are $\partial_{\theta}$ derivatives handled the same way for all three schemes? If so, perhaps we could say this here in a few sentences and/or reference a paper that discusses what is done. ]} \gk{[Added below]}
Note that the WENO discretization is only applied in the 
radial ($\rho$) direction, not the angular ($\theta$) direction. This is 
because the solution is expected to have a very smooth profile in the $\theta$ direction. In that direction, standard centered-differencing 
is used for all methods under consideration in this work.

\subsubsection{Sixth Order Finite Difference Method}
Standard sixth-order centered finite-difference stencils and an eighth-order Kreiss-Oliger dissipation operator are explicitly 
included below. The sixth-order derivative stencil we use takes this form
\begin{equation}
u'(\rho) = \frac{u_{j+3} -9\,{u_{j+2}} + 45\,{u_{j+1}} - 45\,{u_{j-1}} + 9\,u_{j-2}  -\,u_{j-3} }{60\,\Delta \rho}. 
\end{equation}
The Kreiss-Oliger dissipation operator of the proper order for a sixth-order scheme is computed as~\cite{Kreiss73} 
\begin{equation}
Q = \frac{\sigma h^{7} D_{+}^{4}D_{-}^{4}}{256} 
\end{equation}
where $D_{+}=(u_{j+1}-u_{j})/\Delta \rho$ and $D_{-}=(u_{j}-u_{j-1})/\Delta \rho$ are the standard forward and backward differencing operators and 
$\sigma$ is parameter that is usually a value set between $(0,1)$. Based on many numerical experiments, we found that a value of 
$\sigma=0.01$ achieves stability for a large class of computations. 

\subsubsection{WENO Methods}
To discretize the problem \[ u_t + f(u)_\rho =0 \] in space using a WENO method we split the flux into its positive and negative parts
\[ f(u) = f^{+}(u) + f^{-}(u),\]
such that
\[ \frac{d f^{+}(u)}{du}  \geq 0,  \; \; \; \; \mbox{and}  \; \; \; \; \frac{d f^{-}(u)}{du}   \leq 0.\]

\begin{itemize}
\item The third order WENO method: 
For both the positive direction flux $\hat{f}^{+}$ and the negative direction flux  $\hat{f}^{-}$  the smoothness measurements are:
\[ IS_1  =  \left( f^+_{j} - f^+_{j-1} \right)^2,  \; \; \; \; \; \; 
IS_2 =   \left( f^+_{j+1} - f^+_{j} \right)^2
\]
and the weights of the candidate stencils are given by
\[ \alpha_1 = \frac{1}{3} \left( \frac{1}{\epsilon + IS_1} \right)^2, \; \; \; \; \; 
 \alpha_2 = \frac{2}{3} \left( \frac{1}{\epsilon + IS_2} \right)^2,\]
\[ \omega_1 = \frac{\alpha_1}{\alpha_1+\alpha_2}, \; \; \; \; \; \;   \omega_2 = \frac{\alpha_2}{\alpha_1+\alpha_2}.\]

The fluxes are:
\begin{eqnarray*} 
\hat{f}^{+}_{j+ \frac{1}{2}}  &=&  \omega_1 \left( \frac{3}{2} f^+_{j}  - \frac{1}{2} f^+_{j-1} \right)+
\omega_2  \left(  \frac{1}{2} f^+_{j} + \frac{1}{2}  f^+_{j+1} \right) \\
\hat{f}^{-}_{j+ \frac{1}{2}}  &=&  \omega_1 \left( \frac{3}{2} f^+_{j}  - \frac{1}{2} f^+_{j+1} \right)+
\omega_2  \left(  \frac{1}{2} f^+_{j} + \frac{1}{2}  f^+_{j-1} \right) .
\end{eqnarray*}

\item The fifth order WENO method: 
For the positive direction flux $\hat{f}^{+}$, the smoothness measurements are:
\begin{eqnarray*}
IS_0^+ & = & \frac{13}{12} \left( f^+_{j-2} -2 f^+_{j-1} + f^+_{j} \right)^2
          + \frac{1}{4} \left( f^+_{j-2} -4 f^+_{j-1} + 3 f^+_{j} \right)^2 \\
IS_1^+ & = & \frac{13}{12} \left( f^+_{j-1} -2 f^+_{j} + f^+_{j+1} \right)^2
          + \frac{1}{4} \left( f^+_{j-1} - f^+_{j+1} \right)^2 \\
IS_2^+ & = & \frac{13}{12} \left( f^+_{j} -2 f^+_{j+1} + f^+_{j+2} \right)^2
          + \frac{1}{4} \left( 3 f^+_{j} -4 f^+_{j+1} +  f^+_{j+2} \right)^2 
\end{eqnarray*}
For the negative-direction flux,  $\hat{f}^{-}$  the smoothness measurements are:
 \begin{eqnarray*}
IS_0^- & = & \frac{13}{12} \left( f^-_{j+1} -2 f^-_{j+2} + f^-_{j+3} \right)^2
          + \frac{1}{4} \left( 3 f^-_{j+1} -4 f^-_{j+2} +  f^-_{j+3} \right)^2\\
IS_1^- & = & \frac{13}{12} \left( f^-_{j} -2 f^-_{j+1} + f^-_{j+2} \right)^2
          + \frac{1}{4} \left( f^-_{j} - f^-_{j+2} \right)^2 \\
IS_2^- & = & \frac{13}{12} \left( f^-_{j-1} -2 f^-_{j} + f^-_{j+1} \right)^2
          + \frac{1}{4} \left( f^-_{j-1} -4 f^-_{j} + 3 f^-_{j+1} \right)^2 
\end{eqnarray*}

Next, we  calculate the weights for each of the candidate stencils in $\rho$:
$$ \alpha^\pm_0 = \frac{1}{10} \left( \frac{1}{\epsilon + IS_0^\pm}\right)^2 \; \; \; \; \;
\alpha^\pm_1 = \frac{6}{10} \left( \frac{1}{\epsilon + IS_1^\pm}\right)^2 \; \; \; \; \;
\alpha^\pm_2 = \frac{3}{10} \left( \frac{1}{\epsilon + IS_2^\pm}\right)^2  $$
and
\[\omega^\pm_0 = \frac{\alpha^\pm_0}{\alpha^\pm_0 + \alpha^\pm_1 + \alpha^\pm_2} \; \; \; \; \;
\omega^\pm_1 = \frac{\alpha^\pm_1}{\alpha^\pm_0 + \alpha^\pm_1 + \alpha^\pm_2} \; \; \; \; \;
\omega^\pm_2 = \frac{\alpha^\pm_2}{\alpha^\pm_0 + \alpha^\pm_1 + \alpha^\pm_2}.
\]

Finally, we compute the fluxes, comprised of the candidate stencils with their weights:
\begin{eqnarray*} \hat{f}^{+}_{j+ \frac{1}{2}}  &=& \omega^+_0 
\left(\frac{2}{6} f^+_{j-2} -\frac{7}{6} f^+_{j-1} + \frac{11}{6} f^+_{j} \right)
+ \omega^+_1 
\left( - \frac{1}{6} f^+_{j-1} + \frac{5}{6} f^+_{j} + \frac{2}{6} f^+_{j+1} \right)\\
&+& \omega^+_2
\left(\frac{2}{6} f^+_{j} + \frac{5}{6} f^+_{j+1} - \frac{1}{6} f^+_{j+2} \right) 
\end{eqnarray*}
and
\begin{eqnarray*}
 \hat{f}^{-}_{j+ \frac{1}{2}} & = &\omega^-_2 
\left( - \frac{1}{6} f^-_{j-1} + \frac{5}{6} f^-_{j} + \frac{2}{6} f^-_{j+1} \right)
+ \omega^-_1 
\left(\frac{2}{6} f^-_{j} + \frac{5}{6} f^-_{j+1} - \frac{1}{6} f^-_{j+2} \right) \\
&+& \omega^+_0
\left(\frac{11}{6} f^-_{j+1} - \frac{7}{6} f^-_{j+2} + \frac{2}{6} f^-_{j+3} \right).
\end{eqnarray*}
\end{itemize}
Now we discretize each flux as follows:
\[ f^{+}(u)_\rho = \frac{1}{\Delta \rho} \left(
\hat{f}^{+}_{j+ \frac{1}{2}} - \hat{f}^{+}_{j- \frac{1}{2}} \right), \; \; 
\mbox{and} \; \; \;
f^{-}(u)_\rho = \frac{1}{\Delta \rho} \left(
\hat{f}^{-}_{j+ \frac{1}{2}} - \hat{f}^{-}_{j- \frac{1}{2}} \right) .\]

For the problem we are solving,  we need to evaluate the  term $ \mbox{\boldmath{$M$}} \partial_{\rho}\mbox{\boldmath{$u$}} $. 
The structure of $\mbox{\boldmath{$M$}} \partial_{\rho}\mbox{\boldmath{$u$}} $ suggests that the first two elements in 
$\mbox{\boldmath{$u$}} $ have wavespeed $b$, while the second two elements have wavespeed $-b$. This is not completely correct, as
the terms $m_{31}$ and $m_{32}$ have an impact as well; however,  these terms are very small compared
to the wavespeed for most of the computational domain. For simplicity, then, we simply
approximate each element of $\partial_{\rho}\mbox{\boldmath{$u$}} $ with a WENO discretization corresponding to the sign of the wavespeed. This convenient shortcut makes the numerical computation significantly faster  and easier to code; however, the argument above is not rigorous and we do not have a complete numerical analysis that  guarantees convergence of this process.

Proceeding with this  approach, we begin by computing the smoothness measurements for each stencil.
Note that each of the four elements of $\mbox{\boldmath{$u$}} $ is in fact a two-dimensional array, with the $\rho$ and $\theta$ directions.
For convenience, we do not explicitly denote the $\theta$ direction.

The determination of  whether to use an upwind or downwind flux depends on whether the corresponding diagonal term in $ \mbox{\boldmath{$M$}}$
is $b$ or $-b$. We treat the first two elements of  $\mbox{\boldmath{$u$}} $ (which correspond to the value of $b >0$
in the diagonal of  $ \mbox{\boldmath{$M$}} $) as we would treat $f^{+}(u)$, and the last two elements of $\mbox{\boldmath{$u$}} $,
(which correspond to the value of $-b <0$ in the diagonal of  $ \mbox{\boldmath{$M$}} $)  as we would treat $f^{-}(u)$.
 For  the first two elements (arrays) of $\mbox{\boldmath{$u$}} $, we define each of them in turn as an array
$\hat{f}^{+}$, where the values  $\hat{f}^{+}_j$ refer to the $\rho$ values. The WENO derivative is approximated for each value of $\theta$, but the 
index refers to the grid-value of $\rho$. 
For  the last two elements of $\mbox{\boldmath{$u$}} $, we define each of them in turn as
 $\hat{f}^{-}$.  The WENO process gives us the arrays $ \hat{f}^{\pm}_{j+ \frac{1}{2}}$ for all gridpoints $\rho_j$, and then we compute the 
 first two arrays of  $\partial_{\rho}\mbox{\boldmath{$u$}} $ by
\[ \frac{1}{\Delta \rho} \left( \hat{f}^{+}_{j+ \frac{1}{2}} - \hat{f}^{+}_{j- \frac{1}{2}} \right) \]
and the last two element of  $\partial_{\rho}\mbox{\boldmath{$u$}} $ are computed by
\[ \frac{1}{\Delta \rho} \left( \hat{f}^{-}_{j+ \frac{1}{2}} - \hat{f}^{-}_{j- \frac{1}{2}} \right). \]

%{\bf Double check direction and that I am not mixing up the directions!}

\subsubsection{A mixed precision implementation of WENO methods}
The WENO algorithm can be seen as having two parts: 
the inexpensive computation of a finite difference approximation to the derivative on several candidate stencils, 
and a costly nonlinear computation of the stencil weights. 
The computational bottleneck is in the computations of the stencil weights used to combine the stencil-based approximations. 
However, while the differentiation requires highly accurate computation, the calculation of the stencil weights does not require high precision. 
As long as the weights add up to $1.0$ in high precision, the weights may not need to be high precision in the smooth regions, especially
if the region where high precision is needed is near the location of the horizon
%\scott{[Q: does refer to the singularity in the Kerr metric? Or near the particle? 
%Also, if the weights need to add up to 1 in high precision (which is quad?) is this possible if each individual weights are computed in double? ]} 
%\gk{[Fixed; should have been horizon. Yes,imagine computing w1, w2 in double .. promoting them to quad .. and computing w3 = 1-w1-w2 in quad.]}.
With this in mind we modified our quad-precision code to carry out a double-precision computation of the weights, then promoting
them to quad-precision before finally assembling the WENO fluxes. This strategy is expected to speed up the computation significantly.

\subsection{Time-discretization}
When WENO is used to semi-discretize a problem of the form
\[ u_t + f(u)_\rho = 0  \]
we obtain a system of ODEs of the form
\[ u_t = F(u) .\]
The WENO method is designed to have an essentially non-oscillatory property when coupled with the forward Euler method,
\[ u^{n+1}  =  u^n +   \Delta t F\left( u^{n} \right)  \]
under some stability condition $\Delta t \leq \Delta t_{FE}$.
To preserve this property, we use a higher order strong stability preserving time discretization, which can be written
as convex combinations of forward Euler schemes. 
Such time-stepping methods will preserve the properties of the spatial discretization coupled with forward Euler, 
under the modified time-step restriction
\[ \Delta t \leq  {\cal{C}}  \Delta t_{FE}. \]
If ${\cal{C}} > 0 $ we call the method {\em Strong Stability Preserving} (SSP).
While the time-step depends on both the spatial and temporal discretizations, 
we isolate the contribution of the temporal discretization to the time-step restriction by 
considering the ratio  ${\cal{C}}$ of the allowable time-step of the high order method to the forward Euler time-step. 
This ratio is called the strong stability preserving coefficient. 
Using this approach, we view the time-step restriction  as a combination of two factors:
 the forward Euler time-step $\Delta t_{FE}$ that comes from the spatial discretization, 
 and the SSP coefficient ${\cal{C}}$  is a property only of the time-discretization. 
 Among methods of similar types, a more relevant quantity is the effective SSP coefficient
${\cal{C}}_{eff} = {\cal{C}}/s$, which takes into account the computational cost of the method at each iteration, 
defined by the number of stages $s$ (typically also the number of function evaluations per time-step). 
In this work we consider two such methods \cite{gottlieb2011strong}.

\noindent{\bf The three stage, third order strong stability preserving Runge--Kutta method SSP-RK(3,3):} 
%\scott{[What is the $F$ operator? Is this the RHS of the semi-discrete scheme?]} \gk{[Sigal? -- I added an explanation to the first paragraph]}
\begin{eqnarray*}
u^{(1)}  &=&  u^n +   \Delta t F\left( u^{n} \right)  \\ 
u^{(2)}  &=& \frac{3}{4} u^n +  \frac{1}{4}  \left( u^{(1)}  + \Delta t F\left( u^{(1)}  \right)  \right)\\ 
u^{n+1}  &=& \frac{1}{3} u^n + \frac{2}{3}  \left( u^{(2)}  + \Delta t F\left( u^{(2)}  \right)  \right).
\end{eqnarray*}
This method has SSP coefficient ${\cal{C}} =1$ and effective SSP coefficient ${\cal{C}}_{eff} =1/3$.

\noindent{\bf The low storage, ten stage, fourth order  strong stability preserving Runge--Kutta method 
SSP-RK(10,4):}
\begin{eqnarray*}
u^{(1)}  &=&  u^n +   \frac{1}{6} \Delta t F\left( u^{n} \right) , \\ 
u^{(i)}  &=&  u^{(i-1)} +   \frac{1}{6} \Delta t F\left( u^{(i-1)} \right)  , \; \; \;  \mbox{for i=2:4}  \\ 
u^{(5)}  &=&   \frac{3}{5} u^n +  \frac{2}{5} u^{(4)} +   \frac{1}{15}  \Delta t F\left( u^{(4)}  \right),  \\ 
u^{(i)}  &=&  u^{(i-1)} +   \frac{1}{6} \Delta t F\left( u^{(i-1)} \right)  , \; \; \;  \mbox{for i=6:9}  \\ 
u^{n+1}  &=& \frac{1}{25} u^n +  \frac{9}{25} u^{(4)}  + \frac{3}{5} u^{(9)} 
+ \frac{3}{50} \Delta t F\left( u^{(4)}  \right)  +  \frac{1}{10} \Delta t F\left( u^{(9)}  \right) .
\end{eqnarray*}
This method has SSP coefficient ${\cal{C}} =6$ and effective SSP coefficient ${\cal{C}}_{eff} =0.6$.

For the problems considered in this paper, we found no significant difference between the performance of the two.

\section{Results}
Aretakis' rigorous results apply to fairly generic situations, as long as the black hole is extremal. In particular, 
the initial wave that perturbs an otherwise isolated extremal black hole could be fairly generic as long as it has 
support on the horizon. For the following results we chose a narrow Gaussian radial profile i.e., a ``wave-packet'' 
centered at $\rho = 1.0$ with a width of $0.22$. The angular profile of this initial perturbation is chosen to be 
the $(\ell,m)=(2,0)$ spherical harmonic. And we present results for the gravitational perturbation case with $s=-2$. 
%\scott{[For the comparison between the 2 WENO methods and FDM, how many grid points are used? Presumably to compare the methods something (e.g number of grid points) is held constant between the different methods. It might be good to point this out here.]}\gk{[fixed]}
%\scott{[What timestepper is used? We say both ssp-rk methods are comparable, but perhaps for completeness we can say which is used in practice to compile the results of this section.]} \gk{[fixed]}
The numerical simulations presented here used a grid size of $ 16384 (\rho) \times 64 (\theta)$ and the SSP-RK(3,3) time-stepper. 

We wrote the Teukolsky equation for a Kerr black hole with parameters $M,a$ for the variable $\Phi$ %\scott{[This seems to be the first time $\Phi$ appears. Also, I didn't quite get the motivation for $\Phi$ -- is the idea simply that its directly related to the charge? Or does it have an interpretation independent of its appearance in the charge computation?]} 
%\gk{[Yeah .. the basic idea is that it has fall-off similar to the scalar case, so it is the natural quantity to use to extend the notion of Aretakis' scalar charge to the GW case. This is too technical for this paper; we would not get into it.]}, 
which is related to the Teukolsky function $\Psi$ in the Kinnersley tetrad and Boyer-Lindquist coordinates via $\Phi=(r/\,\Delta^2)\,\Psi$, 
where $\Delta=r^2-2Mr+a^2$. One may relate this function to the better-known Weyl curvature in the following manner. 
The Weyl curvature scalar $\psi_4^{\rm HH}$ in the Hartle-Hawking tetrad is related to its Kinnersley tetrad counterpart,  
$\psi_4^{\rm K}$, via a type-III transformation, or  
$\psi_4^{\rm HH}=4(r^2+a^2)^2\,\Delta^{-2}\, \psi_4^{\rm K}$ \cite{poisson2004absorption} and that 
$\Psi =(r-ia\,\cos\theta)^4\,\psi_4^{\rm K}$ \cite{teukolsky1973perturbations} 
we find that  
\begin{equation}
\Phi=\frac{r\left(r-ia\,\cos\theta\right)^4}{4\left(r^2+a^2\right)^2}\,\psi_4^{\rm HH}\, ,
\end{equation}
and use $\Phi$ with $\ell=2,m=0$ and $a=M$. 

Previous work by Aretakis and others~\cite{Aretakis_2013,Angelopoulos:2018yvt,burko2020scalar} show that the radial derivatives of the physical field 
on the horizon would be a conserved quantity, a so-called ``charge''. 
%Such results were derived through rigorous mathematical 
%analysis using advanced techniques in differential geometry and PDE theory. 
In addition, the same analysis also leads to an expectation 
of the time-dependent behavior for the physical field itself on the horizon (proportional to inverse time) and for higher radial 
derivatives (proportional to an increasing positive power of time: second-derivative $\propto \tau$, third-derivative $\propto \tau^2$ 
and so on). To summarize:
\begin{equation}
\Phi^{(p+1)}\propto \tau^p\, ,
\end{equation}
when evaluated at the horizon an extremal black hole. It is this rigorous result that we will compare our numerical solutions with, 
thus offering a precise assessment of the quality of our computational results. We will show that the WENO(5,3) method performs the 
best amongst the different methods we tested.

\begin{figure}[h!]
\centering
\includegraphics[scale=0.75]{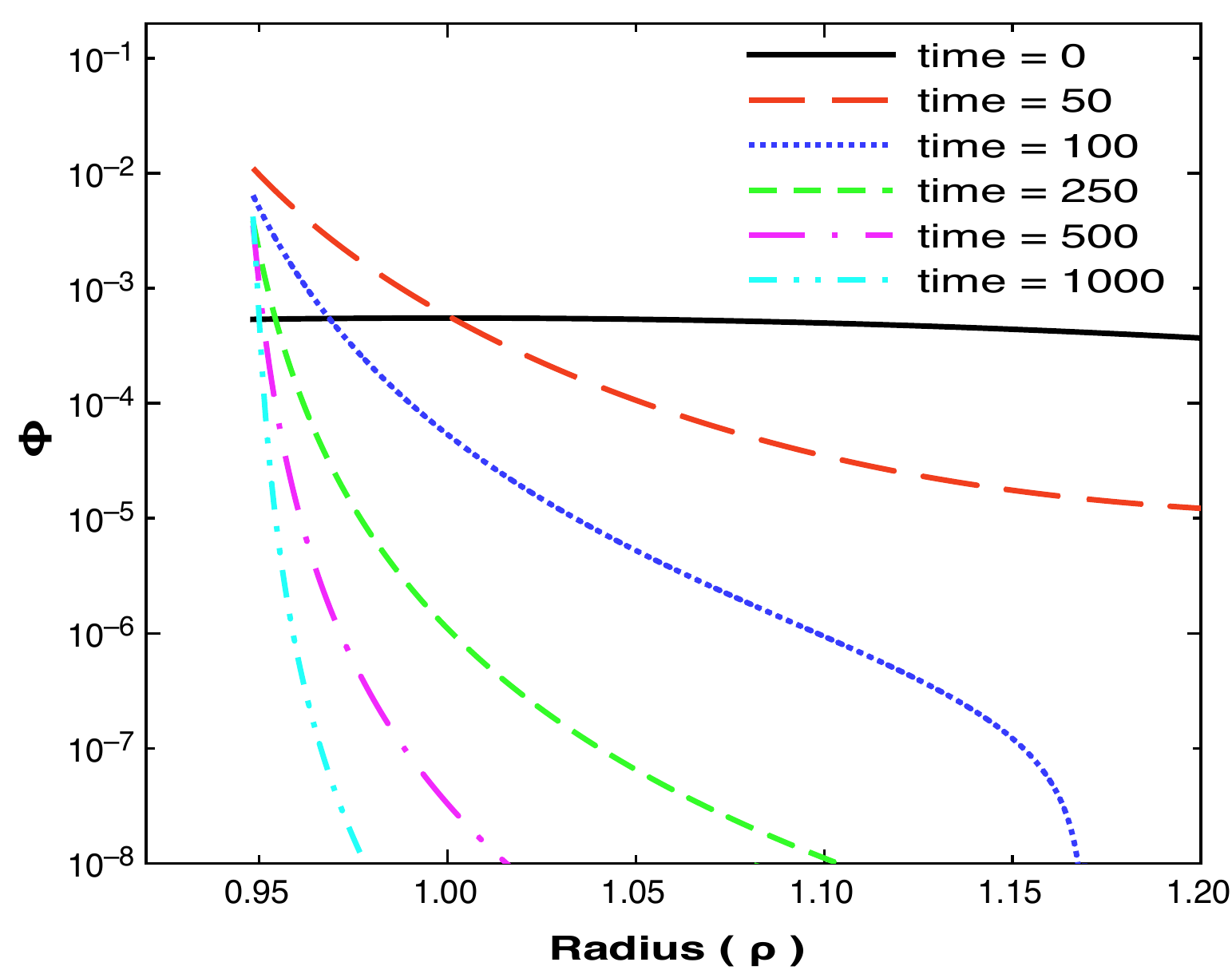}
\caption{The gravitational field spatial profile snapshots in the background of an extremal Kerr black hole: 
The solution's radial dependence at different moments of time. It is clear that as the evolution advances, the 
solution at the horizon (left-end of the computational domain) develops an increasingly sharp profile with the 
higher gradients growing unboundedly.}
\label{fig:radial}
\end{figure}

We begin with an example of some high-quality results from the WENO(5,3) method. We plot $\Phi$  for a fixed $\theta$ as a function 
of $\rho$ for a set of $\tau$ values in Fig.~\ref{fig:radial}. 
These radial snapshots of the solution $\Phi$ at different moments of time suggest how the solution evolves forming 
a sharper-and-sharper feature at the horizon. This plot is highly suggestive of the fact that the higher radial derivatives 
of the field on the horizon will not decay, even though the field itself decays everywhere. 

\begin{figure}[h!]
\centering
\includegraphics[scale=0.75]{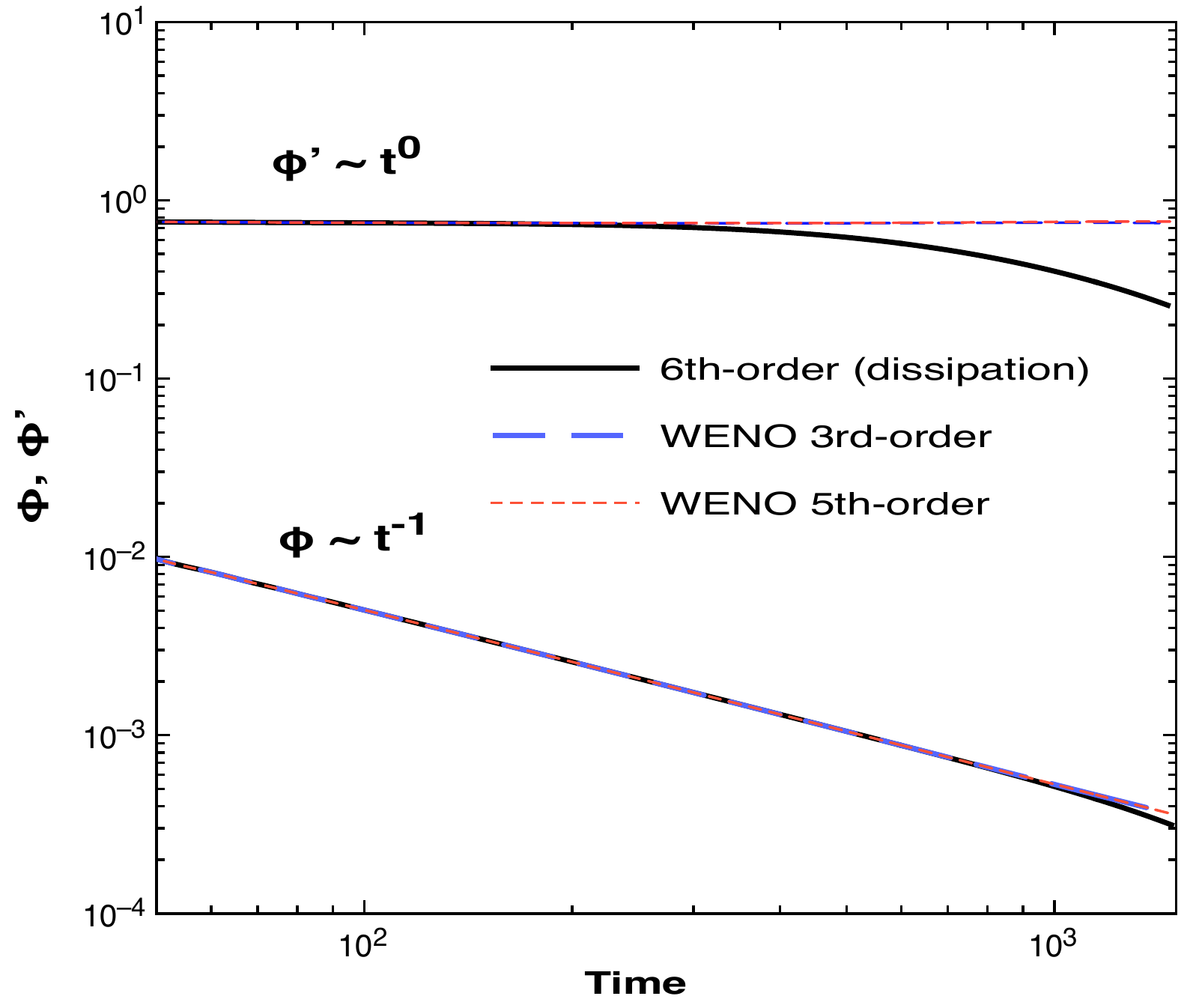}
\caption{The gravitational field $\Phi$, $\Phi'$ temporal profile on the horizon of an extremal Kerr black hole: 
The expected behaviors are indicated in the figure: $\Phi\propto \tau^{-1}$, $\Phi'\propto 1$. It is clear that the 
WENO methods perform significantly better when compared to an even higher-order (dissipative) scheme.}
\label{fig:temporal}
\end{figure}

To observe that more clearly, in Fig.~\ref{fig:temporal} we show the time-dependence of the numerical solution $\Phi$, and its 
radial-derivative $\Phi'$ together. As pointed out above, 
%rigorous mathematical analysis shows that
the physical field  $\Phi\propto \tau^{-1}$
while the radial-derivative $\Phi'$ should be a constant. 
The plot clearly shows that both WENO(3,1) and WENO(5,3) perform significantly better even when compared to a higher-order 
standard (6th-order) finite-difference scheme. The 6th-order standard scheme requires the addition of a small amount of 
dissipation (a standard Kreiss-Oliger 8th-order filter operator~\cite{Kreiss73}) to suppress high-frequency instabilities. 
However, even a small amount of dissipation is enough to significantly degrade the quality of the numerical solution. 

\begin{figure}[h!]
\centering
\includegraphics[scale=0.75]{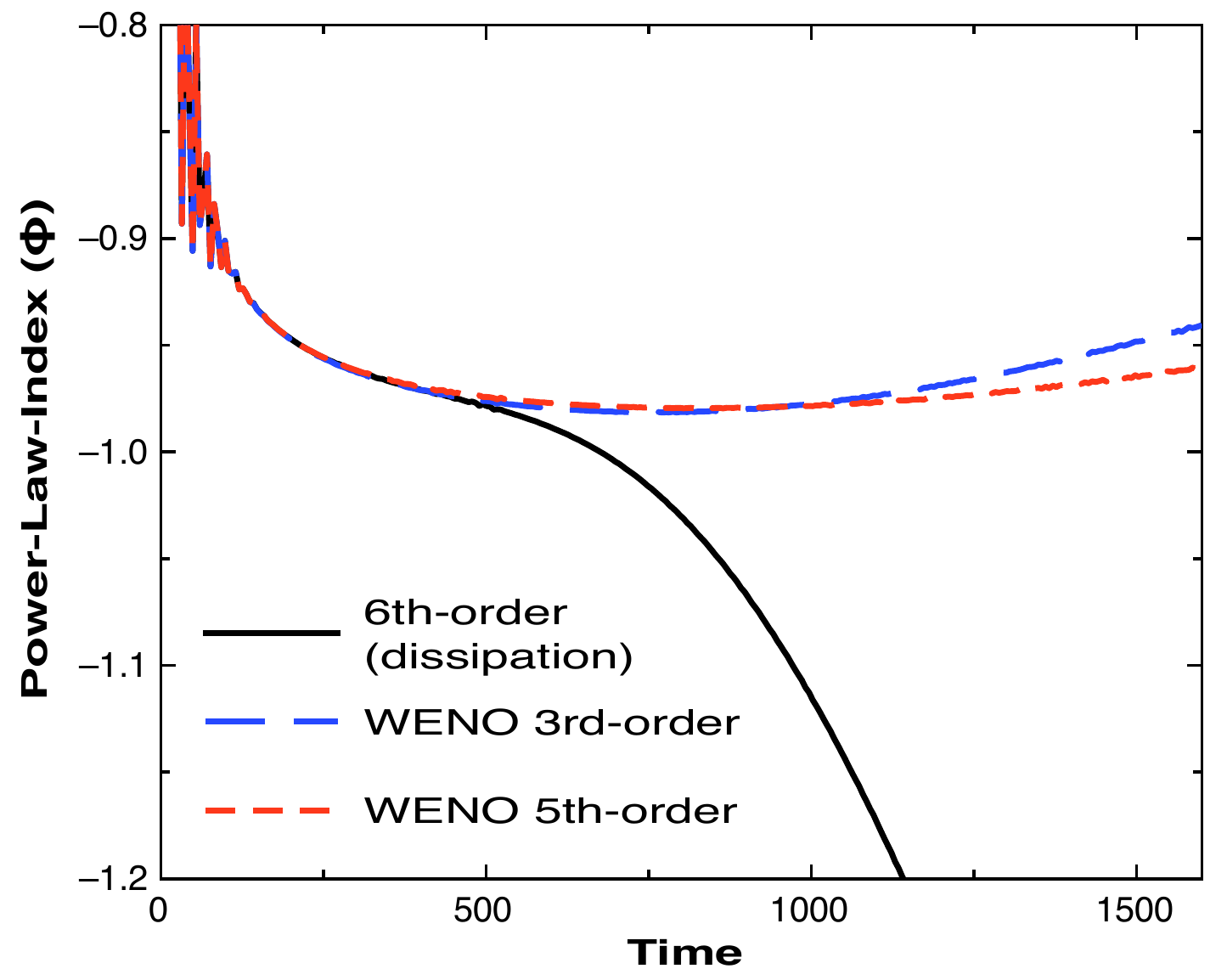}
\caption{The gravitational field $\Phi$ power law tail on the horizon of an extremal Kerr black hole: The expected 
behavior is $\Phi\propto \tau^{-1}$ i.e. a power-law of $-1$. It is clear that the WENO methods perform significantly 
better when compared to an even higher-order (dissipative) scheme. WENO(5,3) performs better than WENO(3,1). }
\label{fig:power-law-field}
\end{figure}

\begin{figure}[h!]
\centering
\includegraphics[scale=0.75]{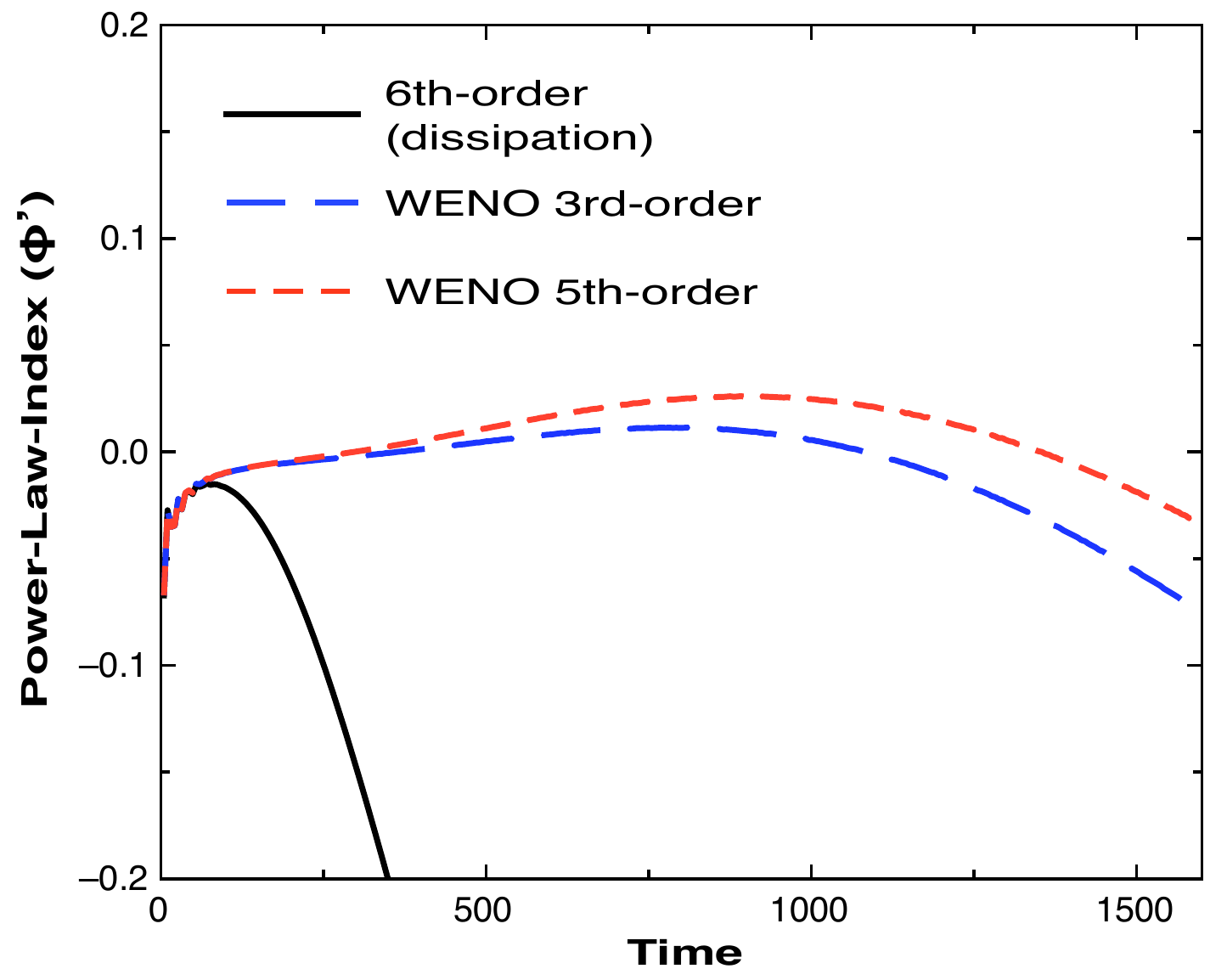}
\caption{The gravitational field $\Phi'$ power law tail on the horizon of an extremal Kerr black hole: The expected 
behavior is $\Phi'\propto \tau^{0}$ i.e. a power-law of $0$. It is clear that the WENO methods perform significantly 
better when compared to an even higher-order (dissipative) scheme. WENO(5,3) performs better than WENO(3,1).}
\label{fig:power-law-der}
\end{figure}

Next, in Fig.~\ref{fig:power-law-field} we show the numerically computed {\em local} power-law index
calculated according to the formula $p = \tau \dot \Phi / \Phi$.
Here we can clearly see that the WENO(5,3) performs the best i.e. maintains a value of $p$ close to $-1$ better 
than WENO(3,1) and also the standard 6th-order finite-difference scheme with dissipation. The same can be seen in 
the $\Phi'$ data as well, as shown in Fig.~\ref{fig:power-law-der} and for higher-order derivatives.

\begin{figure}[h!]
\centering
\includegraphics[scale=0.75]{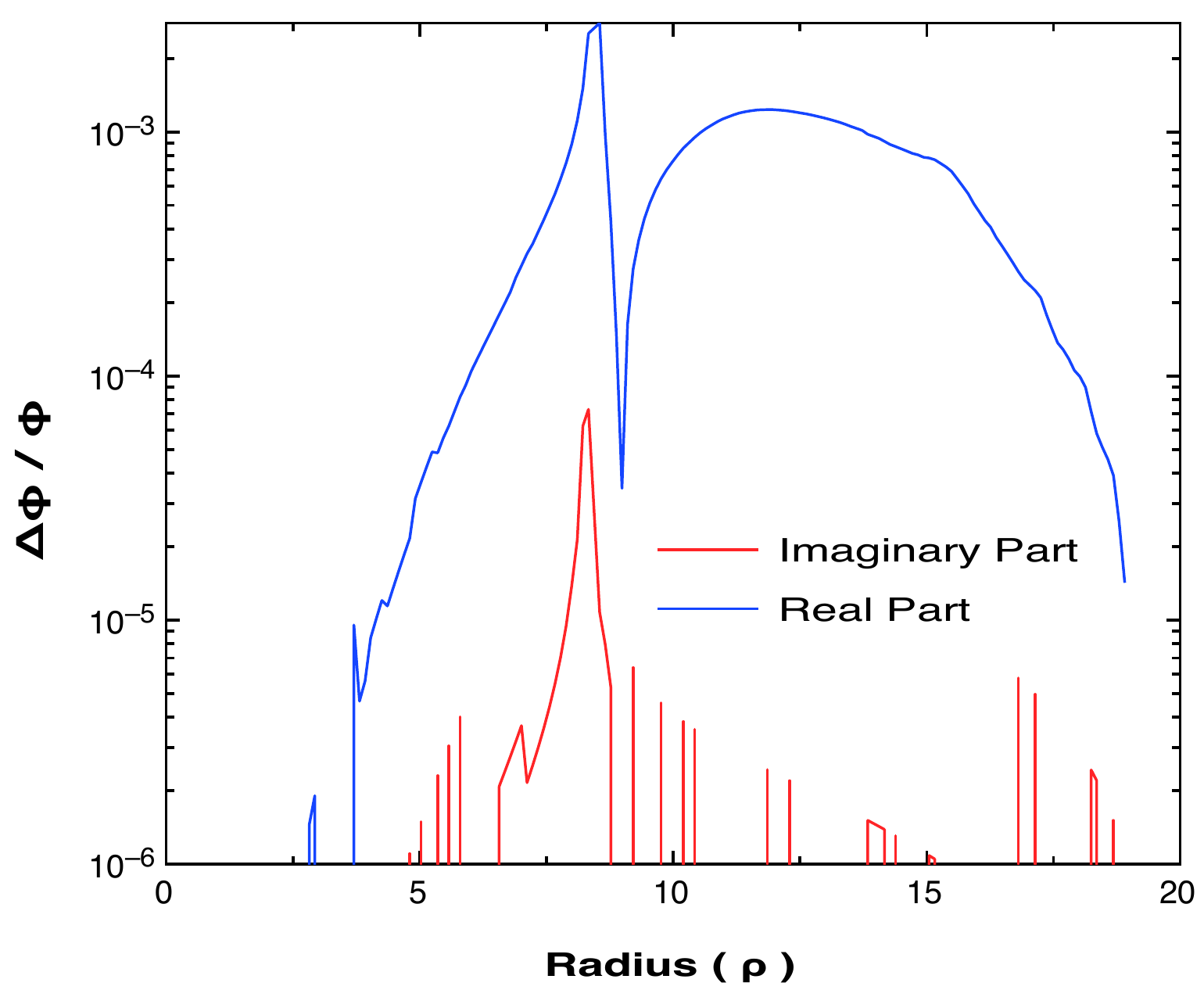}
\caption{Relative differences between two numerical solutions: spatial profile snapshot of the relative differences 
at the {\em last} time-step i.e $\tau=2000$ between two WENO solutions. One solution was generated using {\em full} quadruple-precision 
floating-point computations; the other using a {\em mixed} precision approach i.e. the WENO ``weights'' are computed in 
double-precision. While the relative error becomes significant in the middle portion of the computational domain; it stays 
very low near the horizon while offering a $3.3$ fold speedup.
%\scott{[Q: would it be worthwhile to use a y-axis scale that goes as low as necessary to show the relative error curve's behavior near radius of about 1? Also, is it known what causes the sharp spike around a radius of 9, or why the relative errors reach to such large values in the interiors (naively I would have thought the error would be on the order of double precision). ]}\gk{[Sadly the data in the files was saved just to 6-7 sig. figures. On your other question .. the power-law decay is faster away from the horizon; and as you know, tails computations are numerically very delicate.]}
}
\label{fig:mixed}
\end{figure}

In Fig.~\ref{fig:mixed} we depict the relative difference between two WENO solutions -- one using {\em full} quadruple-precision 
numerics and the other, using a {\em mixed} precision approach i.e. only the WENO ``weights'' are computed using  
double-precision floating-point operations. The mixed precision approach offers a $3.3$ fold speedup while having no impact 
on the results in the regime of interest i.e. near the horizon. Recall that the computed solution $\Phi$ is complex valued. 
It is interesting to note that even though we begin the evolution using purely real initial data, the system's evolution 
introduces a physical phase shift due to the spin of the black hole resulting in a non-zero imaginary part. The mixed precision 
related error in the imaginary part is considerably lower than in the real part. 

\subsection{Parallel Scaling Results}
In this section we present the parallel scaling performance results of our WENO(5,3) code on a GPGPU-cluster with multiple Nvidia V100 GPGPUs. 

First, we offer some details on both the GPGPU {\em many}-core processor architecture and implementation details relevant 
to our code. Nvidia's CUDA framework\footnote{\url{https://developer.nvidia.com/about-cuda}} is a set of software layers that 
allow for GPU devices to become more accessible to the average computational scientist. At a high level, there is a CUDA Runtime 
API; while at a low level, there is the CUDA Driver API. Each function call of the Runtime API is broken into simpler instructions 
and managed by the Driver API. Through CUDA, the GPU (called {\em device}) is accessible to the CPU (called {\em host}) as a 
co-processor with its own memory. The device executes a function (usually referred to as a {\em kernel}) in a ``data parallel'' 
model, which means that a number of threads run the same program on different data.

A data-parallel model is straightforward to implement in a code like ours. We simply perform a domain-decomposition of our 
finite-difference numerical grid and allocate the different parts of the grid to different GPU cores. Each thread computes a 
time-step for a single pair of $\rho$ and $\theta$ grid values. Note that all these calculations are independent, i.e. no  
communication is necessary between the GPU threads.
%\scott{[is some communication necessary for the timesteps or RK stages? Or does each thread get a bit more data than it needs to avoid communication overhead? ]}
%\gk{[On the GPU .. it is a shared memory situation; each thread has ALL the grid data.]}
However, it is necessary to establish the appropriate data communication 
between the GPU-cores and the remaining code that is executing on the CPU respectively.
%\scott{[If its relevant, perhaps it could be good to mention what kinds of communication is needed as the evolution proceeds and/or between stencils. I've somehow had the impression that higher-order WENO is harder to parallelize because the stencil sizes get wider and so there's a lot of data transfer. If the data is availble, it would be interesting to have third and fifth order WENO scaling plots on figure 6. If its not available, perhaps we should mention here which WENO variant is being used.]}
%\gk{[It is 5th order WENO. I think the reason why scaling is good despite the wider stencil is that this is a flops bound computation (quad precision). You're correct that under some
%conditions scaling gets much worse for higher-order methods; something I've argued with Sigal now-and-then. Even wrote a paper on it to convince her  http://worldcomp-proceedings.com/proc/p2014/CSC7132.pdf LOL.]}
%On the Tesla CUDA GPU, we use {\em cudaMemcpy} 
%instructions to achieve the same. In total, only a rather small amount of data is required to be transferred back and forth from 
%the GPU at every time-step of the numerical evolution. Because of the rather modest amount of data-transfer involved, we do not 
%make use of any advanced memory management features on the GPU architecture. 
Of course, this fine-grain decomposition on the GPU is performed after a standard coarse-grain level domain-decomposition using 
MPI over the many GPUs that are part of a tightly-coupled cluster. 

It is worth pointing out that this CUDA implementation of our WENO code is fairly straightforward. It should also 
be mentioned that we do not attempt to hand-tune the codes to tailor them for each architecture, in order to obtain maximal 
performance. Instead, we rely on the mature compiler suites to perform all low-level optimizations (such as vectorization) 
automatically. 

\begin{figure}[h!]
\centering
\includegraphics[scale=0.75]{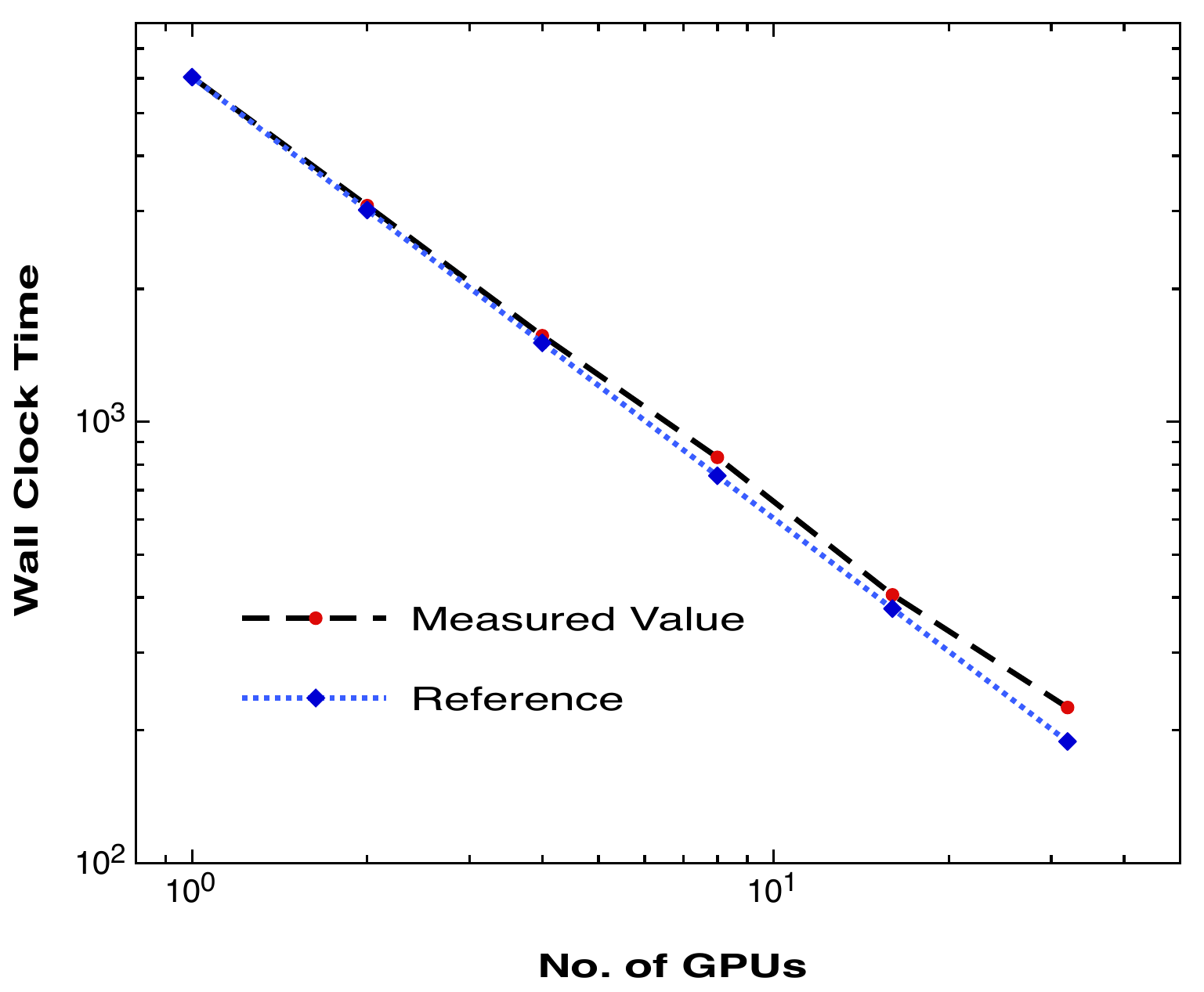}
\caption{Parallel scaling performance of the fifth order WENO method on a GPU cluster. We report the 
total wall-clock time for the same computation using different 
numbers of GPU resources (up to 32). The parallel framework used by the code is CUDA+MPI. The simulations were performed on MIT's 
{\em Satori} supercomputer that has an IBM Power9 CPU + quad Nvidia V100 GPUs node architecture. For reference, perfect scaling 
data is indicated on the same graph. It is clear that the code exhibits near-ideal scaling performance.
We suspect the good scaling performance is partly due to the quad-precision operations,
which in turn increases the computational load of our algorithm relative to the data movement cost expected for wider stencils.
%\scott{[Gaurav -- I've added your explanation here since I think its really helpful and informative for the reader. Please see if this is OK.]}
}
\label{fig:scaling}
\end{figure}

The performance results on MIT's {\em Satori} GPGPU supercomputer are presented in Fig.~\ref{fig:scaling}. This system is an 
IBM machine with 64-nodes; each node offers two Power9 CPUs (32-cores each) and four (4) Nvidia V100 GPGPUs connected via 
NVLink2. We ran the exact same computation using different numbers of GPUs: 1, 2, 4, 8, 16, 32 each time. Note that the 32-GPU 
computation used 8-nodes of this system i.e. a 12.5\% of the entire system size. The results are clearly indicative of nearly 
perfect scaling performance.
%\scott{[For GPU novices like myself, is there some way to quantify the benefits of going from a CPU to a GPU? I'm not sure what the relevant comparison would be -- maybe the time on 1 CPU node vs 1 GPU node? ]} \gk{[done; see below]}
As demonstrated by one of the authors in Ref.~\cite{khanna2013high} the GPU-acceleration alone speeds up the computation by 
nearly a factor of $50$ over a CPU-only computation. 

\section{Conclusion}
In this work we applied a number of GPU-accelerated numerical methods, including a novel mixed-precision fifth order WENO method to solve the Teukolsky equation -- the master equation of black hole perturbation theory. Our results are in three distinct areas: numerical analysis, efficient computing, and gravitational wave astrophysics.

First, we showed that the WENO finite-difference methods of Jiang and Shu \cite{jiang1996efficient} out-performed the sixth order centered difference method with a standard Kreiss-Oliger filter. The dissipation added to the sixth order centered difference is needed to stabilize the method; however, it adversely impacts the quality of the solution so that the expected behavior is not observed near the horizon. The WENO methods correct this issue, as they achieve stability by clever automated stencil choosing and no additional dissipation. For this reason, they attain the correct behavior near the horizon. As expected, the low-order WENO is not sufficiently accurate for this problem, but it still out-performs the sixth order centered difference method with the filter. The fifth order (third order near discontinuities) WENO method we employ confirms the predicted behavior of the system and proves to be appropriate for solving the problem of interest i.e., computing the so-called Aretakis charge, which is an important conserved physical quantity. 

Next, we considered two computational approaches to speed up this code. Regardless of the algorithms used, this is a computationally intensive problem, which requires a very high precision computation (quad precision) to accurately resolve the long-time integration. This problem requires GPU-acceleration to complete in reasonable time. We show that our GPU-accelerated code scales very well as more GPU devices are used. Additionally, to further speed up the computation, we considered a novel mixed-precision approach to the WENO algorithm, by computing the WENO weights in double precision while the rest of the code is computed in quad precision. The mixed precision WENO approach results in a speed-up of a factor of 3.3. We investigated the errors introduced by this approach, and found that the relative error in the solution is very small ($< 10^{-4}$) near the horizon, which is the regime of interest. This minor loss in accuracy is a small and worthwhile price to pay for the dramatic speedup in the code.

%Finally, the scientific results that such high-accuracy simulations enable may result in a significant discovery in the field of computational gravitational wave physics -- the gravitational wave analog of the (scalar) Aretakis charge may be a measurable quantity by observatories like LIGO/Virgo and if indeed detected would prove the existence of extremal or near-extremal black holes. This is a long-standing open problem in the field of black hole physics. 

% slight rewording of paragraph shown above
Finally, such high-accuracy simulations may enable significant discoveries in the field of computational and observational gravitational-wave physics. For example, the gravitational wave analog of the (scalar) Aretakis charge is in principle a measurable quantity that future observatories, like LIGO/Virgo, could detect. Observational evidence for the Aretakis charge would reveal the existence of extremal or near-extremal black holes, settling a long-standing open problem in the field of black hole physics. 

% since this material is already summarized in the conclusion, perhaps it can be removed?
%In conclusion, in this work we compared known numerical methods and showed that the fifth order WENO method is most suitable for important computations in the field of black hole and gravitational wave physics. We presented a GPU-accelerated version of this algorithm, which was further accelerated by a mixed-precision modification to the fifth order WENO method. This combination allowed us to efficiently and accurately compute the behavior of important physical quantities that may ultimately result in a conclusive discovery about the existence of extremal black holes in Nature.

\section{Acknowledgements}
Many of the computations were performed on the MIT/IBM {\em Satori} GPU supercomputer supported by the Massachusetts Green High Performance Computing Center (MGHPCC). 
%{\bf Place holder for Summit}
The authors acknowledge support of NSF Grants No. PHY-2010685 (G.K) 
and No. DMS-1912716 (S.F, S.G, and G.K), AFOSR Grant No. FA9550-18-1-0383 (S.G) and Office of Naval Research/Defense University Research Instrumentation Program (ONR/DURIP) Grant No. N00014181255. 
This material is based upon work supported by the National Science Foundation under Grant No. DMS-1439786 while a subset of the authors were in residence at the Institute for Computational and Experimental Research in Mathematics in Providence, RI, during the Advances in Computational Relativity program. A part of this research is sponsored by the Office of Advanced Scientific Computing Research; US Department of Energy, and was performed at the Oak Ridge National Laboratory, which is managed by UT-Battelle, LLC under Contract no. De-AC05-00OR22725. This manuscript has been authored by UT-Battelle, LLC, under contract DE-AC05- 00OR22725 with the US Department of Energy. The United States Government retains and the publisher, by accepting the article for publication, acknowledges that the United States Government retains a non-exclusive, paid-up, irrevocable, world- wide license to publish or reproduce the published form of this manuscript, or allow others to do so, for United States Government purposes.

\section{Conflict of Interest}
On behalf of all authors, the corresponding author states that there is no conflict of interest. 

\bibliographystyle{unsrt}
\bibliography{references}
%\begin{thebibliography}{11}

%\bibitem{ko_diss} H. Kreiss and J. Oliger, ``Methods for the Approximate Solution of Time Dependent Problems'', 
%Global Atmospheric Research Programme: GARP Publication Series (1973).

%\end{thebibliography}

\end{document}